\documentclass[12pt]{amsart}
\usepackage{amssymb}
\usepackage{amsfonts}
\usepackage{amsmath}
\usepackage{MnSymbol}

\usepackage{stmaryrd}

\newtheorem{example}{Example}[section]

\newtheorem{theorem}[example]{Theorem}

\newtheorem{proposition}[example]{Proposition}

\def\llb{\llbracket}
\def\rrb{\rrbracket}

\def\sym{{\it Sym}}

\def\<{\langle}
\def\>{\rangle}

\def\C{\operatorname{\mathbb C}}

\def\Z{\operatorname{\mathbb Z}}

\def\N{\operatorname{\mathbb N}}
\def\tr{\operatorname{tr}}

\def\SG{{\mathfrak S}}

\def\dim{{\rm dim}}

\def\ch{\operatorname {ch}}

\def\maj{{\rm maj}}

\newdimen\Squaresize \Squaresize=18pt
\newdimen\Thickness \Thickness=0.5pt

\def\Square#1{\hbox{\vrule width \Thickness
   \vbox to \Squaresize{\hrule height \Thickness\vss
      \hbox to \Squaresize{\hss#1\hss}
   \vss\hrule height\Thickness}
\unskip\vrule width \Thickness}
\kern-\Thickness}

\def\Vsquare#1{\vbox{\Square{$#1$}}\kern-\Thickness}

\def\thisbox#1{\kern-.09ex\fbox{#1}}
\def\downbox#1{\lower1.200em\hbox{#1}}

\title[]%
{Stability properties of inner plethyms\\(Lecture Notes)}
\author[J.-Y.~Thibon]{ Jean-Yves Thibon}

\address[]{Laboratoire d'informatique Gaspard-Monge\\
Université Gustave Eiffel, 
5, Boulevard Descartes \\ Champs-sur-Marne \\
77454 Marne-la-Vall\'ee cedex 2 \\
France}
\email[]{jean-yves.thibon@univ-eiffel.fr} 
\date{\today}

\date{\today}

\begin{document}

\begin{abstract}
	The inner plethysm of symmetric functions corresponds to the $\lambda$-ring
	operations of the representation ring $R(\SG_n)$ of the symmetric group.
	It is known since the work of Littlewood that this operation possesses
	stability properties w.r.t. $n$. These properties have been explained in
	terms of vertex operators [Scharf and Thibon, Adv. Math. 104 (1994), 30-58].
	Another approach [Orellana and Zabrocki, Adv. Math. 390
	(2021), \# 107943], based on an expression of character values as symmetric
	functions of the eigenvalues of permutation matrices, has been proposed recently.
	This note develops the theory from scratch,  discusses the link between both
	approaches and provides new proofs of some recent results.
\end{abstract}

\maketitle

\section{Introduction}

The term {\em inner plethysm}, introduced by D.E. Littlewood \cite{Lit1}, refers
to the operation on symmetric functions corresponding to the composition of representations
of the symmetric group $\SG_n$  with representations of the general linear group. For example,
given a linear representation $\rho$ of $\SG_n$ on a vector space $V$, that is, a group
homomorphism $\rho:\ \SG_n\rightarrow GL(V)$, one may consider the representations
$\Lambda^k(\rho)$ in the exterior powers $\Lambda^k(V)$. These operations endow
the representation ring $R(\SG_n)$ with the structure of a $\lambda$-ring \cite{Knu} and
since $R(\SG_n)$ can be identified whith the homogeneous component $Sym_n$ of degree $n$ of the
ring of symmetric funtions $Sym$, this space is itself endowed with a $\lambda$-ring structure,
different from the standard one of $Sym$, induced by the composition of representations of the
general linear groups. This last composition, denoted by $f\circ g$ or $f[g]$, is the usual
(or outer) plethysm, so that for example, the character of $GL(V)$ on the $j$-th exterior power of the
$i$-th exterior power $\Lambda^j(\Lambda^i(V))$ is the plethysm of elementary symmetric functions
$e_j\circ e_i$. Thus, it makes sense to denote the Frobenius characteristic  of the $i$-th exterior 
power of a representation of $\SG_n$ of characteristic $f$ by $\hat e_i[f]$ (inner plethysm of
$f$ by $e_i$), so that its $j$-th exterior power would be $\hat e_j[\hat e_i[f]] =\widehat{e_j\circ e_i}[f]$.

Remarkably, $R(\SG_n)$ is generated  as a $\lambda$-ring by a single element, which can be taken
as the $n$-dimensional vector representation of $\SG_n$ (by permutation matrices) or as its
unique  non-trivial irreducible component, which is of dimension $n-1$. This important result,
which has been rediscovered many times (see, {\it e.g.}, \cite{Mar}), seems to have been first noticed by
P. H. Butler \cite{But}. It implies in particular that any character of the symmetric group can
be expressed as a symmetric function of the eigenvalues of permutation matrices. Such expressions
have been recently investigated by Orellana and Zabrocki \cite{OZ},
Assaf and Speyer \cite{AS} and Ryba \cite{Ry}.
Such expressions imply stability properties, which can also be derived by different methods, such
as vertex operators. 
These notes, which 
correspond roughly to a few talks given over the years at the Combinatorics Seminar in Marne-la-Vallée, will discuss the relations between
the different points of view, and sometimes provide new proofs of old or recent results.

\subsection*{Acknowledgements} This research has been partially supported by the
program CARPLO (ANR-20-CE40-0007) of the Agence Nationale de la Recheche.

\section{Notations and background}

We shall assume that the reader is familiar with the notation of Macdonald's book \cite{Mcd}.

Representations and conjugacy classes of $\SG_n$ are indexed by partitions $\mu$
of $n$, represented as nonincreasing sequences $\mu_1\ge\mu_2\ge\ldots\ge\mu_r>0$
or in exponential notation $(1^{m1}2^{m_2}\cdots n^{m_n})$.

The irreducible representation of $\SG_n$ indexed by $\lambda$ is denoted by $[\lambda]$, and
its character by $\chi^\lambda$. 
Its Frobenius characteristic is the Schur function $s_\lambda=\ch(\chi^\lambda)$, which was written
$\{\lambda\}$ by Littlewood. Littlewood used the symbol $\otimes$ (now reserved for
tensor products) to denote outer plethysms: $\{\lambda\}\otimes\{\mu\}$ is now 
denoted by $s_\mu\circ s_\lambda$ or $s_\mu[s_\lambda]$ (or even $s_\mu(s_\lambda)$).
Littlewood's notation for $\hat s_\lambda[s_\mu]$ was $\{\mu\}\odot\{\lambda\}$.

The pointwise product of central functions translates as the internal product of $Sym$, denoted
by a $*$. On the power-sum basis,
\begin{equation}
	p_\lambda*p_\mu = z_\lambda\delta_{\lambda\mu}p_\lambda.
\end{equation}

Recall that for any $f\in\sym$, there is a differential operator $D_f$
on $\sym$ (the {\it Foulkes derivative\footnote{Denoted by $f^\perp$ in \cite{Mcd}.}}) defined as the adjoint of
the multiplication operator $g\mapsto fg$, that is
\begin{equation}
\< fg\, ,\, h\> = \<g\, ,\,D_f h\> \ .
\end{equation}

Introducing the series 
\begin{equation}
\sigma_z(X):=\sum_{r\ge 0} z^r h_r(X) \quad {\rm and}\quad
\lambda_z(X):=\sum_{r\ge 0} z^r e_r(X) \ ,
\end{equation}
the generating series for the Schur functions indexed by vectors
of the form $(r,\lambda)$ (where $\lambda$ is a fixed partition) can
then be expressed as
\begin{equation}\label{SVERTEX}
\sum_{r\in\Z}z^r s_{(r,\lambda)} = \sigma_zD_{\lambda_{-1/z}} s_\lambda \ .
\end{equation}
This identity is established by expanding the Jacobi-Trudi determinant by
its first row, which causes the appearance of the skew Schur functions
$s_{\lambda/(1^k)}=D_{e_k}s_\lambda$. The operator 
\begin{equation}
\Gamma_z = \sigma_z\, D_{\lambda_{-1/z}}
\end{equation}
is a typical example of the so-called {\it vertex operators}  
(see {\it e.g.} Kac's book \cite{Kac} Chap. 14 for 
examples and references).

In terms of power-sums,
\begin{equation}
\Gamma_z = \exp\left\{\sum_{k\ge 1}\frac{z^k}{k}p_k\right\}\exp\left\{-\sum_{l\ge 1}z^{-l}\frac{\partial}{\partial p_l}\right\}.
\end{equation}

Remark that in $\lambda$-ring notation, 
$$
D_{\sigma_z}f\, (X)=f(X+z)\quad {\rm and}\quad
D_{\lambda_{-z}}f\, (X) =f(X-z)
$$
when $z$ is an element of rank one (which means that
$e_r(z)=0$ for $r>1$), so that
$$
\Gamma_zf\, (X) = \sigma_z(X)\, f\left( X-{1\over z}\right) \ .
$$

Thus, for a partition $\mu$ of $n$, 
\begin{equation}
\begin{split}
\chi^{n-k,\lambda}_\mu &= \<\Gamma_1 s_\lambda,p_\mu\>\\
& = \<\sigma_1D_{\lambda_{-1}}s_\lambda,p_\mu\>\\
&= \left\<\sum_k(-1)^ks_{\lambda/1^k},\prod_{i\ge 1}(1+p_i)^{m_i}\right\>\\
&= \sum_k (-1)^k \sum_{\nu\subseteq\mu}\prod_i{m_i\choose n_i}\<s_{\lambda/1^k},p_\nu\>\\
	&=\Xi^\lambda(m_1,m_2,\cdots,m_k),
\end{split}
\end{equation}
a polynomial in the $m_i$, independent of $n$, which is moreover a $\Z$-linear combination of products
of binomial coefficients ${m_i\choose n_i}$. These have been called {\it character polynomials}
by Specht \cite{Spe}.

As a consequence, there exist stable ($n$-independent) formulas for the reduction of
Kronecker products or inner plethysms. These formulas are stated in terms of the
{\it reduced notation} $\<\lambda\> = [n-|\lambda|,\lambda]$ of Littlewood.
We see that we can identify  $\<\lambda\>$
with the generating series $\Gamma_1 s_\lambda$
and denote it alternatively by $\<s_\lambda\>$, interpreting $\<\cdot\>$ as a linear
operator.


The following identity plays a fundamental role in the derivation of stable character formulas: 
\begin{theorem}\cite{Th0}\label{th:Th0}
	Let $(u_\lambda)$ be any homogeneous basis of $Sym$, and let $(v_\lambda)$
	be its adjoint basis. Then, for any symmetric functions $f,g$,
	\begin{equation}\label{eq:redprod}
	(\sigma_1 f)*(\sigma_1g)=\sigma_1\sum_{\alpha,\beta}D_{u_\alpha}f\cdot D_{u_\beta}g\cdot v_\alpha*v_\beta.
\end{equation}
\end{theorem}

\medskip
{\footnotesize
For example, for the tensor powers of the vector representation, this gives by induction
\begin{equation}
(\sigma_1 h_1)^{*m}=\sigma_1T_m(h_1)
\end{equation}
where $T_m$ are the Touchard polynomials. Indeed, this is true for $m=1$,
and
\begin{equation}
	(\sigma_1 h_1)^{*m}= (\sigma_1 h_1)^{*m-1}=(\sigma_1h_1)=\sigma_1\left(T_{m-1}(h_1)h_1 +D_{s_1}T_{m-1}(h_1)D_{s_1}(h_1)s_1*s_1\right)
\end{equation}
so that $T_m$ satisfies $T_m(x)=x(T_m(x)+T_m'(x))$.
}

\medskip
We set $\llangle f\rrangle=\sigma_1f$. The $\llangle h_\mu\rrangle$, of for short $\llangle\mu\rrangle$,  are called {\it stable
permutation characters} \cite{STW}.

\medskip
{\footnotesize
Thus, the above example reads $\llangle 1\rrangle^{*m}=\llangle T_m(h_1)\rrangle$.
}

\medskip
Theorem \ref{th:Th0} implies the existence of coefficients $\bar g^{\nu}_{\lambda\mu}$
and $\bar d^{\nu}_{\lambda\mu}$ such that
\begin{equation}
	\<\lambda\>*\<\mu\>=\sum_\nu \bar g^{\nu}_{\lambda\mu}\<\nu\>,\
	\llangle\lambda\rrangle*\llangle\mu\rrangle=
\sum_\nu \bar d^{\nu}_{\lambda\mu}\llangle\nu\rrangle,
\end{equation}
called {\it reduced Kronecker coefficients} (see {\it e.g.}, \cite{STW}).

As a consequence, the internal product is well-defined on series of the form $\sigma_1f$, where $f$
is a symmetric function of finite degree. The linear span of these series will be called the ring of stable characters,
and denoted by $\widehat{Sym}$.

\section{Inner plethysm: first steps}
 
Let $V=\C^n$ and $\rho:\ \SG_n\rightarrow GL(V)$ be the representation
by permutation matrices. Its character $\chi(\tau)=\tr \rho(\tau)$
is the number of fixed points of $\tau$:
 if the cycle type of $\tau$
is $\mu=(1^{m_1}2^{m_2}\cdots n^{m_n})$, then 
\begin{equation}
\chi(\tau) = m_1.
\end{equation}
Since $m_1=\<h_1,p_\mu(X+1)\>$, recalling that $D_{\sigma_1}f(X)=f(X+1)$,
we have $m_1=\<h_{n-1,1},p_\mu\>$
so  that  its Frobenius characteristic is $h_{n-1,1}$.


Any symmetric function can be expressed as a polynomial (with rational coefficients)
in the power sums $p_k$. The corresponding operators on  representations are usually called Adams
operations, and denoted by $\psi^k$: for a representation $\pi$ of a group $G$ of character $\xi$, one defines
\begin{equation}
\psi^k(\xi)(g) = \tr \pi(g^k).
\end{equation}
This is in general only a virtual character.
In the case of the vector representation of $\SG_n$, $\psi^k(\chi)(\tau)$ is the number of
fixed points of $\tau^k$. Thus, for $\tau$ of type $\mu$,
\begin{equation}\label{eq:psik}
\psi^k(\chi)(\tau) = \sum_{d|k}dm_d.
\end{equation}
The Frobenius characteristic of this virtual character is thus $\hat p_k[h_{n-1,1}]$.

Note that $dm_d=\<p_d,p_\mu(X+1)\>$,
so that \eqref{eq:psik} is equivalent to
\begin{equation}
\sum_{n\ge 1}	\hat p_k[h_{n-1,1}]=\sigma_1\sum_{d|k}p_d.
\end{equation}
Thus, all the $m_i$, hence also the character polynomials, can be expressed as
inner plethysms of $m_1$. This already proves that, as a $\psi$-ring, $R(\SG_n)$ is
generated by the vector representation $V$.

\medskip

{\footnotesize
The first examples are
\begin{align}
\<1\> &\leftrightarrow \Xi^1 = \<s_1(X-1),p_\mu(X+1)\> = m_1-1 \\
\<2\> &\leftrightarrow \Xi^2 = \<s_2(X-1),p_\mu(X+1)\> =m_2 + {m_1\choose 2}-m_1\\
\<11\>&\leftrightarrow \Xi^{11} =\<s_{11}(X-1),p_\mu(X+1)\>={m_1\choose 2}-m_2-m_1+1 
\end{align}
from which we can compute
\begin{align}
\<1\>*\<1\> &\leftrightarrow (m_1-1)^2 = \Xi^2+\Xi^{11}+\Xi^1+\Xi^0\\
\hat p_2\<1\>&\leftrightarrow \psi^2(m_1-1)=2m_2+m_1-1\\
\hat h_2\<1\>&\leftrightarrow \frac12(\psi^2(m_1-1)+\psi^{11}(m_1-1)) =m_2+{m_1\choose 2}=\Xi^2+\Xi^1+\Xi^0
\end{align}
so that
\begin{align}
\<1\>*\<1\> &= \<2 \> + \<11 \> + \<1 \> + \<0 \>\\
\hat h_2\<1\> &=   \<2 \> +  \<1 \> + \<0 \>\\
\hat e_2\<1\> &= \<11\>
\end{align}
and we can check for example that
\begin{equation}
s_{41}*s_{41}= s_{{5}}+s_{{41}}+s_{{32}}+s_{{311}},\quad \hat e_2(s_{41})= s_{311}.
 \end{equation}
}

\medskip
Next, we can form the generating series
\begin{equation}
	\hat\sigma_x[h_{n-1,1}](\tau)=
	\exp\left\{\sum_{k\ge 1}\frac{x^k}{k}\sum_{d|k}dm_d\right\}
=	\prod_{k\ge 1}(1-x^k)^{-m_k(\tau)},
\end{equation}
after rearranging the sum in the exponential.
This provides the expression of $\hat h_k[h_{n-1,1}](\tau)$ as a polynomial in the $m_i$.
Also,
\begin{equation}
	\hat\lambda_{-x}[h_{n-1,1}](\tau)=\prod_{k\ge 1}(1-x^k)^{m_k(\tau)}.
\end{equation}
This is of course the (reciprocal) characteristic polynomial of the permutation matrix of $\tau$,
and the calculation could have been done the other way round.

In terms of symmetric functions, this remark allows the computation of $\ch \Lambda^k(\rho)=  \hat e_k[h_{n-1,1}]$.
Indeed, the (reciprocal) characteristic polynomial of $\rho(\tau)$ is
\begin{equation}
|I-x\rho(\tau)| = \sum_{k=0}^n (-x)^k \tr \Lambda^k(\rho(\tau)),
\end{equation}
and since the reciprocal characteristic polynomial of a $p$-cycle is $1-x^p$,
if follows from the cycle decomposition of $\tau$ that
\begin{equation}
|I-x\rho(\tau)| = \prod_i (1-x^i)^{m_i} = p_\mu[1-x].
\end{equation}
The Frobenius characteristic of $\Lambda^k(\rho)$ is therefore the coefficient
of $(-x)^k$ in
\begin{equation}
\sum_{k=0}^n (-x)^k \hat e_k[h_{n-1,1}] = \sum_{\mu\vdash n}p_\mu[1-x]\frac{p_\mu}{z_\mu} = h_n[(1-x)X].
\end{equation}
Now, 
\begin{equation}
h_n[(1-x)X] = \sum_{k=0}^n h_{n-k}(X)h_{k}(-xX) = \sum_{k=0}^n(-x)^kh_{n-k}e_k
\end{equation}
whence 
\begin{equation}
 \hat e_k[h_{n-1,1}]  = h_{n-k}e_k.
\end{equation}
The nontrivial irreducible component of $V$ is $[n-1,1]$, and
writing $s_{n-1,1}=h_{n-1,1}-h_n$, we have as well
\begin{equation}
	\hat\lambda_{-x}[h_{n-1,1}-h_n]=\hat\lambda_{-x}[h_{n-1,1}]*\hat\sigma_x[h_n] =\frac{\hat\lambda_{-x}[h_{n-1,1}]}{1-x}=\frac{h_n[(1-x)X]}{1-x}
\end{equation}
since $\hat\sigma_x[h_n]=(1-x)^{-1}h_n$,
and taking into account the well-known expansion
\begin{equation}
h_n[(1-x)X] = \sum_{k=0}^n(1-x)(-x)^ks_{n-k,1^k},
\end{equation}
we arrive at
\begin{equation}
\hat e_k[s_{n-1,1}]=s_{n-k,1^k}.
\end{equation}

Define the alphabet $\Omega_\mu$ by the condition
\begin{equation}
\lambda_{-x}(\Omega_\mu) = p_\mu(1-x)
\end{equation}
{\it i.e.}, $\Omega_\mu$ is the multiset consisting of the eigenvalues of a permutation
of type $\mu$, so that the trace of such a permutation on $\Lambda^kV$ is $e_k(\Omega_\mu)$.
We have therefore
\begin{equation}
e_k(\Omega_\mu) = \<h_{n-k}e_k,p_\mu\> 
	=\<\sigma_1e_k,p_\mu\>=\<e_k,D_{\sigma_1}p_\mu\>
	=\<e_k,p_\mu(X+1)\>,
\end{equation}  

\section{The representation of $\SG_n$ on polynomials}

The traces of the symmetric powers being the coefficients of the inverse
of the (reciprocal) characteristic polynomial, we have
\begin{equation}
\sum_{k\ge 0} x^k \tr S^k(\rho)(\tau) = p_\mu\left[\frac{1}{1-x}\right]
\end{equation} 
so that the graded characteristic of $S(V)=\C[x_1,\ldots,x_n]$ is
\begin{equation}
\ch_x \C[x_1,\ldots,x_n] = h_n\left[\frac{X}{1-x}\right].
\end{equation}
Its expansion on Schur functions is known only through its expansion
on ribbon skew Schur functions\footnote{For a composition $I=(i_1,\ldots,i_r)$, 
$$
r_I =\left|\begin{matrix}h_{i_1}&h_{i_1+i_2}&h_{i_1+i_2+i_3}&\cdots & h_{i_1+\cdots+i_r}\\
	                   1    &    h_{i_2}&h_{i_2+i_3}   & \cdots & h_{i_2+\cdots i_r}\\
			   0    &       1   &  h_{i_3}     &  \cdots& h_{i_3+\cdots +i_r}\\
			   \vdots&     &       & \ddots        &  \vdots\\
0    &0           &0            &\cdots    &h_{i_r}\end{matrix}\right|.
$$}
\begin{equation}\label{eq:rib}
 h_n\left[\frac{X}{1-x}\right] = \frac1{(x)_n}\sum_{I\vDash n}x^{\maj(I)} r_I
	=\frac1{(x)_n}\sum_{\lambda\vdash n}f^\lambda(x)s_\lambda
\end{equation} 
where $f^\lambda(x)$ is the generating function by major index of the standard tableaux of shape $\lambda$.
As the orbit of each monomial spans a permutation representation, it is also
interesting to write down the expansion on the basis $h_\mu$. Its generating series is
\begin{equation}
	\sigma_t \left[\frac{X}{1-x}\right] =\sum_{n\ge 0}t^n h_n\left[\frac{X}{1-x}\right]
	=\prod_{k\ge_0}\sigma_{tx^k}(X)
\end{equation}
so that
\begin{equation}
	\hat h_k[h_{n-1,1}]=\sum_{m_0+m_1+\cdots+m_k=k,\atop m_1+2m_2+\cdots+nm_n=n=n}h_{m_0}h_{m_1}\cdots h_{m_n}.
\end{equation}

The characteristic of the orbit of a monomial $x^\mu$ is $h_{m_0}h_{m_1}\cdots h_{m_n}$,
where $m_0+m_1+\cdots+m_n=n$.

{\footnotesize
\begin{example}{\rm
For $n=3$, 
	\begin{align*}
		\mu=(100)&\rightarrow h_{21}\\
		\mu=(200)&\rightarrow h_{21}\\
		\mu=(110)&\rightarrow h_{21}\\
		\mu=(300)&\rightarrow h_{21}\\
		\mu=(210)&\rightarrow h_{111}\\
		\mu=(111)&\rightarrow h_{3}
	\end{align*}
	so that 
	$$\hat h_1[h_{21}] = h_{21},\ \hat h_2[h_{21}] = 2h_{21},\ 
	\hat h_3[h_{21}] = h_3+h_{21}+h_{111}.$$ 
	}
\end{example}
}

In terms of stable characters,
\begin{equation}
	\hat h_k[\sigma_1h_1]=\sigma_1\sum_{m_1+2m_2+\cdots+km_k=k}h_{m_1,m_2,\ldots,m_k}.
\end{equation}

\section{Littlewood duality}

Let $X^{(i)}=\{x^{(i)}_1,x^{(i)}_2,\ldots x^{(i)}_n\}$ be $r$ sets of variables,
and consider the tensor product
\begin{equation}
W = \C[X^{(1)}]\otimes_{\C\SG_n}\C[X^{(2)}]\otimes_{\C\SG_n}\cdots\otimes_{\C\SG_n}\C[X^{(r)}]
\end{equation}
Since the graded characteristic of a single polynomial ring $\C[X]$ is $h_n\left(\frac{X}{1-q}\right)$,
the $r$-graded characteristic of $W$ is
	\begin{multline}\label{eq:multgr}
		h_n\left(\frac{X}{1-q_1}\right)*h_n\left(\frac{X}{1-q_2}\right)*\cdots * h_n\left(\frac{X}{1-q_r}\right)\\
		=h_n\left(\frac{X}{(1-q_1)(1-q_2)\cdots (1-q_r)}\right)= h_n[\sigma_1(Q)X],
	\end{multline}
where $Q=\{q_1,q_2,\ldots,q_r\}$.
This is the term of degree $n$ in
\begin{equation}
\sigma_1[\sigma_1(Q)X] = \sum_{\alpha\in\N^r}q^\alpha \hat h_\alpha[\sigma_1 h_1] = \sum_{\mu}m_\mu(Q)\hat h_\mu[\sigma_1 h_1],
\end{equation}
and taking a scalar product of this expression with any $g\in\sym$, we have
\begin{equation}
\begin{split}
	\sum_{\mu}m_\mu(Q)\<\hat h_\mu[\sigma_1 h_1],g\>&=\<\sigma_1[\sigma_1(Q)X),g(X)\>\\
	&=g[\sigma_1(Q)]=\sum_\mu\<h_\mu,g[\sigma_1]\>m_\mu(Q),
\end{split}
\end{equation}
so that
\begin{equation}
\<\hat h_\mu[\sigma_1 h_1],g\> = \<h_\mu,g[\sigma_1]\>.
\end{equation}

By linearity, we obtain the following statement, relating inner
and outer plethysms. 

\begin{theorem}\label{th:lit}
For any
two symmetric functions $f,g$,
\begin{equation}
\<\hat f[\sigma_1h_1],g\> = \<f,g[\sigma_1]\>.
\end{equation}
\end{theorem}
This is Littlewood's duality \footnote{Other proofs can be found in \cite{ST} and \cite{Stan}} (Theorem XI of \cite{Lit1}).

\medskip
{\footnotesize
Note that combining \eqref{eq:multgr} with \eqref{eq:rib}, 
and taking into account the relation of the internal product
to the descent algebra,
we obtain the multigraded
Hilbert series of the invariants (multisymmetric functions) as
\begin{equation}
\sum_{\alpha\in\N^r}{\bf q}^\alpha\dim Sym_\alpha^{n,r}=
\<\sigma_1[\sigma_1(Q)X),h_n(X)\>
=\sum_{\sigma_1\circ\cdots\circ\sigma_r =id;\ \sigma_i\in\SG_n}
\frac{q_1^{\maj \sigma_1}\cdots q_r^{\maj \sigma_r}}{(q_1)_n\cdots (q_r)_n}
\end{equation}
where the sum runs over all $r$-factorisations of the identity in $\SG_n$
(A. M. Garsia and I. Gessel, Advances in Math. {\rm 31}
(1979),  288--305).
}

\section{Weight spaces}

Theorem \ref{th:lit} describes in particular the branching rule $GL(n,\C)\downarrow \SG_n$,
where $\SG_n$ is embedded as the subgroup of permutation matrices: the multiplicity
of the irreducible representation $[\mu]$ of $\SG_n$ in the restriction of the
irreducible representation $V_\lambda$ of $GL(n,\C)$ is equal to
$\<s_\lambda, s_\mu[\sigma_1]\>$.

When $\lambda\vdash kn$, the weight space $V_\lambda(k,k,\ldots,k)$ is stable
under the action of $\SG_n$, and its characteristic can be computed by a formula
of Gay \cite{Gay} which is somewhat similar to Theorem \ref{th:lit}. Under restriction
to $SL(n)$, this is  the zero weight space.

To derive it, le us rather start from a product of symmetric powers
\begin{equation}
S^\lambda(V):=S^{\lambda_1}(V)\otimes S^{\lambda_2}(V)\otimes \cdots\otimes S^{\lambda_r}(V)
\end{equation}
whose $GL(n)$-character is $h_\lambda$. The elements of this space can be interpreted as polynomials
in $r$ sets of $n$ variables $X^{(i)}$ as above, which are homogeneous of degree $\lambda_i$ 
for each set $X^{(i)}$. 

The zero weight space is spanned by monomials which are homogeneous of degree $k$ in each set of variables
$X_{i}:=\{x_i^{(j)},\ j=1,\ldots,n\}$, which can be represented by nonnegative integer matrices
with row sums $\lambda$ and column sums $(k^n)$. The symmetric group acts by permuting the columns
of these matrices, hence by a permutation representation. 

Let us say that such a matrix has type $\mu=(1^{m_1}2^{m_2}\cdots n^{m_n})$ if 
it has $m_i$ columns $C_i$, with the $C_i$ distinct. The orbit of such a matrix
is then a permutation representation of characteristic $h_\mu$. 

The possible columns, which must have sum $k$, can be encoded by the monomials
of $h_k(Q)$ over an auxiliary alphabet  $Q=\{q_1,\dots,q_r\}$ as above. The number of
matrices of type $\mu$ is therefore equal to the coefficient of $m_\lambda(Q)$
in $m_\mu[h_k(Q)]$, i.e. to $\<h_\lambda,m_\mu[h_k]\>$. Thus, the 
restriction to the zero weight space is given by the adjoint $F_k^\dagger$
of the linear operator $F_k:\ f\mapsto f[h_k]$:

\begin{proposition}[\cite{Gut, Gay}]If $\lambda$ is a partition of $nk$, the Frobenius characteristic of the action of
$\SG_n$ on the zero weight space of the simple module $V_\lambda$
of $SL(n,\C)$ is given by
\begin{equation}
\<\ch V_\lambda(0)\downarrow\SG_n, s_\mu\>= \<s_\lambda, s_\mu[h_k]\>.
\end{equation}
\end{proposition}

{\footnotesize
For example, the zero weight space of $S^{321}(\C^3)$ is spanned by the orbits of the monomials
corresponding to the matrices
\begin{equation}
\begin{pmatrix}2&0&1\\0&2&0\\0&0&1\end{pmatrix}\rightarrow h_{111},\quad
\begin{pmatrix}2&1&0\\0&1&1\\0&0&1\end{pmatrix}\rightarrow h_{111}\quad\text{and}\quad
\begin{pmatrix}1&1&1\\1&1&0\\0&0&1\end{pmatrix}\rightarrow h_{21}
\end{equation}
so that $\ch S^{321}(\C^3)\downarrow\SG_3=h_{21}+2h_{111}$, and one can check that
$\<h_{321},s_\mu[h_2]\>=1$ for $\mu=(21)$, $=2$ for $\mu=(111)$, and $=0$ for $\mu=3$.

For the irreducible module $V_{321}$, the result is $F_2^\dagger(s_{321})=s_{21}$.
}

\medskip
The other weight spaces are not stable under $\SG_n$, but the direct sum of their orbits are.
Denoting for a module $M$ by $\overline{M(\nu)}$ the direct sum
$\bigoplus_{\alpha\in \SG_n(\nu)}M(\alpha)$, 
the same reasoning shows that the characteristic of the restriction of $S^\lambda(\C^n)$ to 
$\overline{S^\lambda(\C^n)(\nu)}$ is the coefficient of $t^\nu$ in 
\begin{equation}\label{eq:ws}
\sum_{\mu\vdash n}\<h_\lambda, m_\mu[t_0+t_1h_1+t_2h_2+\cdots]\>h_\mu
\end{equation}
and by linearity, for an irreducible representation $V_\lambda$,
\begin{equation}
\sum_{\mu\vdash n}\<s_\lambda, s_\mu[t_0+t_1h_1+t_2h_2+\cdots]\>s_\mu =\sum_{\nu\vdash n}\ch \overline{V_\lambda(\nu)}\downarrow_{\SG_n}t^\nu.
\end{equation}
This is equivalent to \cite{NTV}[Cor. 2]. Such decompositions are obtained in \cite{Har} by first restricting to the subgroup of
monomial matrices.
Note that \eqref{eq:ws} is a common generalization of Littlewood’s duality and of Gay’s
formula.

\medskip
{\footnotesize
For example, the restriction of $V_{321}$ of $GL(3)$ to $\SG_3$ decomposes according to the orbits of the weigths as
$$  
t_{1} t_{2} t_{3}s_{111} + \left(t_{2}^{3} + 2 t_{1} t_{2} t_{3}\right)s_{21} + t_{1} t_{2} t_{3}s_{3}
$$
For $S^{321}(\C^3)$, one finds
$$
\left(2 t_{2}^{3} + 12 t_{1} t_{2} t_{3} + 3 t_{0} t_{3}^{2} + 3 t_{1}^{2} t_{4} + 5 t_{0} t_{2} t_{4} + 3 t_{0} t_{1} t_{5}\right)h_{111} + \left(t_{2}^{3} + 2 t_{1}^{2} t_{4} + t_{0}^{2} t_{6}\right)h_{21}
$$

As another example, let us reproduce Table 1 of \cite{NTV}. Set $t_0=1$. The restrictions
to $\SG_n$ of the orbit spaces of $S^{111}(\C^n)$ are given by the vector partitions 
$$
\begin{pmatrix}1  \\1  \\ 1 \end{pmatrix}\
\begin{pmatrix} 1&0  \\1 &0  \\0 &1  \end{pmatrix}\
\begin{pmatrix} 1&0  \\0 &1  \\1 &0  \end{pmatrix}\
\begin{pmatrix} 1&0  \\ 0&1  \\0 &1  \end{pmatrix}\
\begin{pmatrix} 1&0 &0 \\0 &1 &0 \\0 &0 &1 \end{pmatrix}
$$
so that the restriction to the orbit space $\mu$ is the coefficient of $t_\mu$ in
$$t_3\llangle 1\rrangle+3t_2t_1\llangle 11\rrangle+t_1^3\llangle 111\rrangle.$$
For $S^{21}(\C^n)$, we have the matrices
$$
\begin{pmatrix}2  \\1  \end{pmatrix}\
\begin{pmatrix}2 &0  \\0 &1  \end{pmatrix}\
\begin{pmatrix}1 &1  \\1 &0  \end{pmatrix}\
\begin{pmatrix} 1&1 & 0\\0 &0 &1  \end{pmatrix}
$$
giving 
$$t_3\llangle 1\rrangle+2t_2t_1\llangle 11\rrangle+t_1^3\llangle 21\rrangle.$$
Finally, for $S^3(\C^n)$, we have the partitions
$$
(3)\ (2\ 1)\ (1\ 1\ 1)
$$
giving
$$t_3\llangle 3\rrangle +t_2t_1\llangle 11\rrangle +t_1^3\llangle 3\rrangle.$$

This last representation is irreducible, $S^3=S_3$, so converting the stable permutation
characters into stable characters, we get for $V_3$
$$t_3\<s_3(X+1)\>+t_2t_1\<s_1(X+1)^2\> +t_1^3\<s_3(X+1)\>$$
yielding
$$
t_3(\<1\>+\<0\>) +t_2t_1(\<2\>+\<11\>+2\<1\>+\<0\>)+t_1^3(\<3\>+\<2\>+\<1\>+\<0\>)
$$
which reproduces the first column of \cite[Table 1]{NTV}.

For the second column, we write $s_{21}=h_{21}-h_3$, which gives
$$t_2t_1\llangle 11\rrangle +t_1^3\llangle 21-3\rrangle$$
and in terms of stable characters, this is
$$t_2t_1(\<2\>+\<11\>+2\<1\>+\<0\>)+t_1^3(\<21\>+\<2\>+\<11\>+\<1\>).$$
Finally, writing $s_{111}=h_{111}-3h_{21}+h_3$, we obtain ater the same reductions the last column
in the form
$$
t_1^3(\<111\>+\<11\>).$$
}

\medskip
As observed in \cite{HSW}, if we denote by $b^\lambda_\mu$ the coefficient of $s\mu$
in $\hat s_\lambda[h_{n-1,1}]$, then,
\begin{equation}
\sum_{\lambda\vdash n}b^\lambda_\lambda = \sum_{\lambda\vdash n}\<s_\lambda,s_\lambda[\sigma_1]\>=
\sum_{\lambda\vdash n}\<h_\lambda,m_\lambda[\sigma_1]\>
\end{equation}
is the number of functional patterns (endofunctions) over a set of $n$ elements.
And indeed, this is the dimension of the subspace of $V^{\otimes n}$ of invariants under the 
action of $\SG_n$ given by
\begin{equation}
\sigma:\ e_{i_1}\otimes e_{i_2}\otimes\cdots\otimes e_{i_n}\mapsto
e_{\sigma(i_{\sigma^{-1}(1)})}\otimes e_{\sigma(i_{\sigma^{-1}(2)})}\otimes\cdots\otimes e_{\sigma(i_{\sigma^{-1}(n)})}.
\end{equation}
The orbits of the weight spaces are stable for this action, and endofunctions, regarded as equivalence classes
of words of length $n$ over $[n]$ can be classified according to their weight, that is, the partition formed  
by the number of occurences of each letter.

\medskip
{\footnotesize
 For example, with $n=3$, the stable decompositions abov give
\begin{align}
S^{111}(\C^3)&\rightarrow t_3h_{21}+(3t_2t_1+t_1^3)h_{111}\\
S^{21}(\C^3)&\rightarrow (t_3+t_1^3)h_{21}+2t2t_1h_{111}\\\
S^3(\C^3)&\rightarrow t_3h_{21}+t_2t_1h_{111}+t_3h_3
\end{align}
so that
\begin{equation}
\sum_{\lambda\vdash n}\<h_\lambda,m_\lambda[1+t_1h_1+t_2h_2+\cdots]\>
= 3t_1^3+3t_2t_1+t_3
\end{equation}
corresponding to 3 orbits of weight $(1,1,1)$ (the conjugacy classes of permutations,
represented by the words 123, 132, 231), 3 orbits of weight $(2,1)$ (represented by the
words 112, 122, 121), and one orbit of weight $(3)$ (represented by 111).

For $n=4$, using the alternate expression in terms of Schur functions, one can read from
\cite[Table 2]{NTV} that  decomposition of the next number 19 is
\begin{equation}
5t_1^4+7t_2t_1^2+3t_2^2+3t_3t_1+t_4.
\end{equation}

}

\section{The Butler-Boorman theorem}

It follows in particular from Theorem \ref{th:lit} that the coefficient of 
$h_\nu$ in $\hat h_\mu[\sigma_1 h_1]$ is
\begin{equation}
\<\hat h_\mu[\sigma_1 h_1],m_\nu\> = \<h_\mu,m_\nu[\sigma_1]\>,
\end{equation}
and since
\begin{equation}
m_\nu[\sigma_1] = \sum (X^{\alpha_1})^{\nu_1}(X^{\alpha_2})^{\nu_2}\cdots(X^{\alpha_s})^{\nu_s}
\end{equation}
where the $X^{\alpha_i}$ run over all distinct monomials in $X$, the coefficient of $m_\mu$
in this expression is equal to the coefficient of its leading monomial 
$X^\mu=x_1^{\mu_1}\cdots x_r^{\mu_r}$. Encoding such a monomial by a column vector, itself regarded
as an indeterminate, and ordering these columns  lexicographically
we see that the coefficient of $m_\nu$ is independent
of the first part $\nu_1$ (exponent of the monomial 1), and 
is equal to the number of packed $r\times s$ matrices of nonnegative
integers, up to permutation of the columns (vector partitions), with row sums vector $\mu$, such that the multiplicities
of the different columns are the $\nu_j$ for $j>1$. 

{\footnotesize
\begin{example}\label{ex:h21}\rm
The complete expansion of $\hat h_{21}[\sigma_1h_1]$ is
$\llangle h_{21}+2h_{11}+h_1\rrangle$, which can be read on the matrices
$$
\begin{pmatrix}
1&1&0\\
0&0&1
\end{pmatrix}\
\begin{pmatrix}
1&1\\
1&0
\end{pmatrix}\
\begin{pmatrix}
2&0\\0&1
\end{pmatrix}\
\begin{pmatrix}
2\\1
\end{pmatrix}
$$  
\end{example}
}

From this description, it is clear that 
\begin{equation}
\hat h_\mu[\sigma_1 h_1] = \llangle h_\mu + F_\mu\rrangle,
\end{equation}
where $F_\mu$ is a sum of terms $h_\nu$ with $|\nu|<|\mu|$. This proves
that all stable permutations characters $\llangle h_\mu\rrangle$ can be expressed
as integral linear combinations of inner plethysms $\hat h_\nu[\sigma_1 h_1]$,
hence that 
\begin{theorem}[\cite{But,Boo}]\label{th:bb}
$R(\SG_n)$ is generated as a $\lambda$-ring by $h_{n-1,1}$,
or as well by $s_{n-1,1}$. 
\end{theorem}

But this proves more: these expressions are actually independent of $n$.

\bigskip
These results have been proposed as a method of evaluating inner plethysms in terms
of classical operations on symmetric functions \cite{STW}. To evaluate $\hat f[g]$, first
express $g$ as $g=\hat G[h_{n-1,1}]$. Then, 
\begin{equation}
	\hat f[g] = \hat f[\hat G[h_{n-1,1}]] = \widehat{(f\circ G)}[h_{n-1,1}].
\end{equation}

\section{Reduced notation and eigenvalues of permutation matrices}

Another consequence 
of Theorem \ref{th:bb}
is that any character value $\chi(\tau)$
on a permutation of type $\mu$ is a symmetric function of the
eigenvalues of $\tau$ on $\C^n$, that is, of the alphabet $\Omega_\mu$.
Following \cite{OZ}, define $\tilde s_\lambda$ and $\tilde h_\lambda$ by the conditions
\begin{equation}
	\chi^{n-|\lambda|,\lambda}_\mu = \tilde s_\lambda(\Omega_\mu),\quad \xi^{n-|\lambda|,\lambda}_\mu = \tilde h_\lambda(\Omega_\mu)
\end{equation}
where $\xi^\lambda$ is the permutation character corresponding to $h_\lambda$.
or, equivalently,
\begin{equation}
	\<\lambda\> =\hat{\tilde s}_\lambda[\sigma_1 h_1],\quad \llangle\lambda\rrangle=\hat{\tilde h}_\lambda[\sigma_1h_1].
\end{equation}

{\footnotesize
For example, 
\begin{equation}
\begin{split}
\<22\> &= \sigma_1 s_{22}(X-1)= \sigma_1(s_{22}-s_{21}+s_{11})\\
&= \cdots+s_{\bar 3 22}+s_{\bar 2 22}+ s_{bar 1 22} +s_{022}+s_{122}+s_{222}+s_{322}+\cdots\\
&= 0 +s_{11}+s_{111}+0+0+s_{222}+s_{322}+\cdots
\end{split}
\end{equation}
and with 
\begin{equation}
\tilde s_{22} =  s_{22} - s_{3}- 2s_{21} + 4s_{11} + 2s_{2}-s_{1}  
\end{equation}
one can check that
\begin{equation}
	\tilde s_{22}[h_{11} ] = s_{11},\  \tilde s_{22}[h_{21}] = s_{111},\ \tilde s_{22}[h_{31}] = 0,\ \tilde s_{22}[h_{41} ] = 0,\ \tilde s_{22}[h_{51} ] = s_{222},\ \tilde s_{22}[h_{61} ] = s_{322}. 
\end{equation}
}

Let us calculate a few examples. We already know that 
\begin{equation}
\begin{split}
	\<1^k\>&=\sum_n s_{n-k,1^k}=\sum_n \hat e_k[h_{n-1,1}-h_n]\\
	&=\sum_n\sum_{i+j=k}\hat e_i[h_{n-1,1}]*\hat e_j[-h_n]\\
	&=\sum_n\sum_{i=0}^k(-1)^{k-i}\hat e_i[h_{n-1,1}]
\end{split}
\end{equation}
so that
\begin{equation}
	\tilde s_{1^k}=\sum_{i=0}^k(-1)^{k-i}e_i = e_k[X-1].
\end{equation}
For the complete functions, we have \cite{STW}
\begin{equation}
	\llangle n\rrangle = \hat F_n[\sigma_1h_1],\ \text{with}\ F_n=\sum_{i+2j=n}(-1)^jh_ie_j.
\end{equation}
Indeed, 
\begin{equation}
	\sum_{n\ge 0}q^n\llangle n\rrangle =\sigma_1[(1+q)X]=\sigma_1\left[\frac{1-q^2}{1-q}X\right]
\end{equation}
and
\begin{equation}
	\sigma_1\left[\frac{1-y}{1-x}X\right]=\hat\sigma_x[\sigma_1h_1]*\hat\lambda_{-y}[\sigma_1h_1]
\end{equation}
so that
\begin{equation}
	\llangle n\rrangle = \hat\sigma_q[\sigma_1h_1]*\hat\lambda_{-q^2}[\sigma_1h_1]=\sum_{i,j\ge 0}q^i(-q^2)^j\widehat{h_ie_j}[\sigma_1h_1]
\end{equation}
whence
\begin{equation}
	\tilde h_n= \sum_{i+2j=n}(-1)^jh_ie_j.
\end{equation}

Let us introduce the shorthand 
 notation $\llb f\rrb:=\hat f[\sigma_1 h_1]$, so that $\<\lambda\>=\llb\tilde s_\lambda\rrb$ and
 $\llangle\lambda\rrangle=\llb\tilde h_\lambda\rrb$.

{\footnotesize
\begin{example}{\rm Let us check Eq. (20) of \cite{OZ}.
	As a symmetric function of the eigenvalues, 
	\begin{align*}
		h_{21}&= \llb h_2\rrb *\llb h_1\rrb =\llangle 2+1\rrangle*\llangle 1\rrangle\\
		&=\sigma_1(h_2+h_1)*\sigma_1 h_1\\
		&= \sigma_1[(h_2+h_1)h_1+(h_1+1)h_1]\ \text{by \eqref{eq:redprod}}\\
		&= \sigma_1[h_{21}+2h_{11}+h_1]\\
		&=\sigma_1[s_{21}+s_3+2s_2+2s_{11}+s_1]\\
		&=\<s_{21}(X+1)+s_3(X+1)+2(s_1+1)^2+s_1+1\>\\
		&=\<s_{21}+s_3+4s_2+3s_{11}+7s_1+4\>
	\end{align*}
}
\end{example}
}

\section{Duality}

Define coefficients $c_\lambda^\mu$ by
\begin{equation}
h_\lambda =\sum_\mu c_\lambda^\mu \tilde h_\mu.
\end{equation}
Then,
\begin{equation}
	\hat h_\lambda[\sigma_1h_1]= \sum_\lambda c_\lambda^\mu\hat{\tilde h}_\mu[\sigma_1h_1]=\sigma_1\sum_\lambda c_\lambda^\mu h_\mu
\end{equation}
so that
\begin{align}
	\<\hat h_\lambda[\sigma_1h_1],g\>&=\sum_\mu c_\lambda^\mu\<\sigma_1h_\mu,g\>\\
	&= \sum_\mu c_\lambda^\mu\<h_\mu,g(X+1)\>\\
	&=\<h_\lambda,g[\sigma_1]\>.
\end{align}
Thus, if $g(X+1)=m_\mu$, that is $g(X)=m_\mu(X-1)$, we obtain
\begin{equation}
	c_\lambda^\mu = \<\hat h_\lambda[\sigma_1h_1],g\> = \<h_\lambda,g[\sigma_1]\>
	=\<h_\lambda,m_\mu[\sigma_1-1]\>,
\end{equation}
so that the dual basis of $\tilde h_\mu$ can be identified with $m_\mu[\sigma_1-1]$.

As a consequence, we can see that $c_\lambda^\mu$ is equal to the number of vector partitions
of $\lambda$ whose multiplicities form the partition $\mu$. 

\medskip
{\footnotesize
For example, the matrices of Example \ref{ex:h21}
can  now be read as
\begin{equation}
	h_{21}=\tilde h_{21}+2\tilde h_{11}+\tilde h_1.
\end{equation}
}

Now, if
\begin{equation}
s_\lambda =\sum_\mu a_\lambda^\mu \tilde s_\mu,
\end{equation}
writing
\begin{equation}
	\tilde h_\mu = \sum_\lambda k_{\lambda\mu}\tilde s_\mu,
\end{equation}
we have
\begin{equation}
	\sigma_1 h_\mu =\hat{\tilde h}_\mu[\sigma_1h_1] =  \sum_\lambda k_{\lambda\mu}\hat{\tilde s}_\mu[\sigma_1h_1] =\sum_\lambda k_{\lambda\mu}\sigma_1s_\mu(X-1)
\end{equation}
so that
\begin{equation}
	k_{\lambda\mu}=\<s_\lambda,h_\mu(X+1)\>=\<\sigma_1s_\lambda,h_\mu\>
\end{equation}
whence
\begin{equation}
	a_\lambda^\mu=\<s_\lambda,\sigma_1[\sigma_1-1]s_\mu[\sigma_1-1]\>.
\end{equation}
The dual basis of $\tilde s_\mu$ can therefore be identified with $\sigma_1[\sigma_1-1]s_\mu[\sigma_1-1]$.
Note that this is not of the form $\sigma_1f$ with $f$ of bounded degree, and that the dual of $\widehat{Sym}$
is  spanned by series of the form  $\sigma_1[\sigma_1-1]f[\sigma_1-1]$ where $f$ is of finite degree.

\section{The Assaf-Speyer formula}

Define now coefficients $b_\lambda^\mu$ by
\begin{equation}
\tilde s_\lambda=\sum_\mu b_\lambda^\mu s_\mu.
\end{equation}
By duality,
\begin{equation}
	s_\mu = \sum_\lambda b_\lambda^\mu \tilde s_\lambda^* = \sigma_1[\sigma_1-1]\sum_\lambda b_\lambda^\mu s_\lambda[\sigma_1-1].
\end{equation}
Recall that the Poincar\'e-Birkhoff-Witt theorem is equivalent to the fact that the universal enveloping algebra
$U(L)$ of the free Lie algebra $L$ on a vector space $V$ is isomorphic, as a $GL(V)$-module, to $S(L)$ and also to
$T(V)$. In terms of $GL(V)$-characters, this amounts to the plethystic identity
\begin{equation}\label{eq:pbw}
	\sigma_1\left[\sum_{n\ge 1}\ell_n\right]=\frac1{1-p_1},
\end{equation}
where
\begin{equation}
	\ell_n =\frac1n\sum_{d|n}\mu(d)p_d^{n/d}
\end{equation}
is the character of $GL(V)$ in the homogeneous component $L_n$ of $L$.

An equivalent form is
\begin{equation}\label{eq:invplet}
	\sigma_1\left[-\sum_{n\ge 1}\ell_n(-X)\right]= 1+X
\end{equation}
(this reflects the Koszul duality between the operads $Com$ and $Lie$).

Set for short $S=\sigma_1-1$ and $M=-L(-X)$, so that $S\circ M=M\circ S = p_1=X$.
We can now write
\begin{equation}
	\lambda_{-1}[S]s_\mu = \sum_\lambda b_\lambda^\mu s_\lambda[S]
\end{equation}
and composing by $M$
\begin{equation}
	\lambda_{-1}[S\circ M][M]s_\mu = \sum_\lambda b_\lambda^\mu s_\lambda[S\circ M]
\end{equation}
that is
\begin{equation}
	\lambda_{-1}[X]s_\mu[M] = \sum_\lambda b_\lambda^\mu s_\lambda(X)]
\end{equation}
so that finally,
\begin{align}
	b_\lambda^\mu&=\left\<s_\lambda(X), \lambda_{-1}s_\mu\left[-\sum_{n\ge 1}\ell_n(-X)\right]\right\>\\
	&=\left\<s_\lambda(X-1),s_\mu\left[-\sum_{n\ge 1}\ell_n(-X)\right]\right\>\\
	&=\left\<s_\lambda(-X-1),s_\mu\left[-\sum_{n\ge 1}\ell_n(X)\right]\right\>\\
	&=(-1)^{|\lambda|+|\mu|}\left\<s_{\lambda'}(X+1),s_{\mu'}\left[\sum_{n\ge 1}\ell_n(X)\right]\right\>
\end{align}
which is essentially Eq. (7) of \cite{AS}.

This can be recast as
\begin{equation}
	\tilde s_\lambda(Y)=(-1)^{|\lambda|}\<s_{\lambda'}(X),\sigma_1(X)\lambda_{-1}[\ell(X)Y]\>
\end{equation}
which suggests the existence of a resolution of Specht modules in terms of Schur modules.
Such a resolution is exhibited in \cite{Ry}. If instead of $\tilde s_\lambda$ we choose to compute
$\tilde x_\lambda$, defined by $\llangle s_\lambda\rrangle = \hat{\tilde x}_\lambda[\sigma_1h_1]$,
we can get rid of the parasitic factor $\sigma_1$, so that
\begin{equation}
	\tilde x_\lambda(Y)=(-1)^{|\lambda|}\<s_{\lambda'}(X),\lambda_{-1}[\ell(X)Y]\>.
\end{equation}
At this point, $\ell(X)Y$ can be interpreted as the $G:=GL(V)\times \SG_n$ character\footnote{Here, $Y=h_1(Y)$ stands for the character
of the vector representation of $GL(n,\C)$ which restricts to the permutation representation of $\SG_n$, which means that
in the expansion of this series, the Schur functions  $s_\mu(Y)$ must be interpreted as $\tilde s_\mu$.  }
of the Lie algebra $\mathfrak{g}:=L(V)\otimes \C^n$, where $\SG_n$ acts on $\C^n$ by permutation matrices.
Then, $\lambda_{-1}[\ell(X)Y]$ is the $G$-equivariant Euler characteristic of the Chevalley-Eilenberg
complex of $\mathfrak{g}$. Explicit calculation of $H^i(\mathfrak{g},\C)$ shows that the $\SG_n$-character
of the multiplicity space of $s_{\lambda'}$ is precisely $\tilde x_\lambda$, which provides the sought
resolution.

This can be rewritten as,
\begin{equation}
	\lambda_{-1}[\ell(X)Y]=\lambda_{-1}(X)\sum_{\mu}(-1)^{|\mu|}\tilde s_{\mu'}(Y)s_\mu(X)=\sum_{\mu}(-1)^{|\mu|}\tilde x_{\mu'}(Y)s_\mu(X).
\end{equation}

Finally, if
\begin{equation}
\tilde h_\lambda = \sum_\mu d_\lambda^\mu h_\mu,
\end{equation}
the same reasoning leads to
\begin{equation}
	d_\lambda^\mu = \<h_\lambda, h_\mu[-L(-X)]\>.
\end{equation}

\section{Coproducts of stable characters}

Let
\begin{equation}
\Delta\tilde s_\lambda =\sum_{\mu,\nu}f^{\mu\nu}_\lambda \tilde s_\mu\otimes\tilde s_\nu.
\end{equation}
Since $Sym$ is self dual,
\begin{equation}
f^{\mu\nu}_\lambda =\<\tilde s_\lambda, \tilde s_\mu^*\tilde s_\nu^*\>,
\end{equation}
knowing that $s_\mu^*=(\sigma_1s_\mu)[\bar\sigma]$ (where $\bar\sigma:=\sigma_1-1$),
we can write
\begin{align}
 \tilde s_\mu^*\tilde s_\nu^*&=\sigma_1[\bar\sigma](\sigma_1s_\mu s_\nu)[\bar\sigma]\\
&=\sigma_1[\bar\sigma]\sum_\alpha c^\alpha_{\mu\nu}(\sigma_1 s_\alpha)[\bar\sigma]\\
&=\sigma_1[\bar\sigma]\sum_\alpha c^\alpha_{\mu\nu}\sum_{\lambda/\alpha\in {\rm HS}}s_\lambda[\bar\sigma]\\
&=\sum_\alpha c^\alpha_{\mu\nu}\sum_{\lambda/\alpha\in {\rm HS}}\tilde s_\lambda^*,
\end{align}
(where HS means horizontal strips),
so that \cite[Th. 4.7]{OZ2}
\begin{equation}
f^{\mu\nu}_\lambda =
\sum_{\lambda/\alpha \in{\rm HS}} c^\alpha_{\mu\nu}
\end{equation}

In a similar way, $\tilde h_\lambda^*=m_\lambda[\bar \sigma]$ implies that
$\tilde h_\lambda$ has the same coproduct coefficients as $h_\lambda$.
Also, $\tilde x_\lambda^* =s_\lambda[\bar\sigma]$, which implies that
 $\tilde x_\lambda$ has the same coproduct as $s_\lambda$ \cite[Prop. 4.5 and Cor. 4.6]{OZ2}.

\section{Products of stable characters}

In \cite{OZ3}, Orellana and Zabrocki establish combinatorial formulas for various products
of stable characters. All these formulas are consequences of the one for $\tilde h_\lambda \tilde s_\mu$,
which can easily be derived from Donin's formula for $\<h_\lambda,s_\mu*s_\nu\>=\<h_\lambda*s_\mu,s_\nu\>$ ({\it cf.} \cite{Th0}).
Indeed, if $\lambda$ is of length $r$,
\begin{equation}
	h_\lambda*s_\mu=\mu_r[(h_{\lambda_1} h_{\lambda_2}\cdots h_{\lambda_r})*_r\Delta^rs_\mu]
\end{equation}
where $\mu_r$ denotes $r$-fold multiplication and $\Delta^r$ is the iterated coproduct
valued in $Sym^{\otimes r}$, so that
\begin{equation}
	\<h_\lambda*s_\mu,s_\nu\> = \sum_{I_1,\ldots,I_r\atop J_1,\ldots,J_r}\<s_{I_1} ,s_{J_1}\> \<s_{I_2} ,s_{J_2}\>\cdots  \<s_{I_r} ,s_{J_r}\>
\end{equation}
where the sum runs overs all the decompositions of $\mu$ and $\nu$ in successive skew diagrams
\begin{equation}
	\mu = I_1 I_2\ldots I_r,\quad \nu=J_1 J_2\ldots J_r,\quad |I_k|=|J_k|=\lambda_k.
\end{equation}
The first diagrams $I_1$ and $J_1$ are partitions, so that the first scalar  product $\<s_{I_1} ,s_{J_1}\>$ 
can be only 1 or 0.
For $\lambda_1$ large enough ($\lambda_1>\max(\mu_1,\nu_1)$), the other skew diagrams will be independent of its
value. Thus, there exist universal coefficients suth that
\begin{equation}
	\llangle \lambda\rrangle*\<\mu\> = \sum_\nu l_{\lambda\mu}^\nu\<\nu\>,
\end{equation}
and 
\begin{equation}
l_{\lambda\mu}^\nu = 
\sum_{I_0, I_1,\ldots,I_r\atop J_0, J_1,\ldots,J_r}\<s_{I_1} ,s_{J_1}\> \<s_{I_2} ,s_{J_2}\>\cdots  \<s_{I_r} ,s_{J_r}\>
\end{equation}
where the $I_k,J_k$ are decompositions asz above corresponding to the partitions
$(N-|\mu|,\mu)$, $(N-|\nu|,\nu)$ and $(N-|\lambda|,\lambda)$ with $N$ large enough.

Each scalar product has a simple combinatorial interpretation, from which 
that of the total coefficient can be easily derived.

\section{Appendix: the free Lie algebra and  the pure braid group}
Equations \eqref{eq:pbw} and \eqref{eq:invplet} for $\sigma_1[L(X)]$ and $\sigma_1[-L(-X)]$ raise the question
of the interpretation of $\sigma_t[L(X)]$ and $\sigma_t[-L(-X)]$. It turns out that 
Equation \eqref{eq:invplet} is  related to the cohomology of the pure braid group $P_n$.
Its homogeneous component of degree $n$ is its equivariant Poincaré characteristic, which is indeed
0 except in the trivial cases $n=0,1$. More interesting is the equivariant Poincaré polynomial
\begin{equation}\label{eq:ppol}
	\sum_{i\ge 0}(-t)^i\ch H^{n-i}(P_n;\C) = \left. \sigma_t[-L(-X)]\right|_{{\rm degree}\ n}.
\end{equation}
The right-hand side can be expanded as 
\begin{equation}
	\begin{split}
		\sigma_t[-L(-X)]& = \prod_{i\ge 1}(1+p_i)^{\ell_i(t)}\\
		&=(1+p_1)^t(1+p_2)^{\frac12(t^2-t)}(1+p_3)^{\frac13(t^3-t)}(1+p_4)^{\frac14(t^4-t^2)}\cdots
	\end{split}
\end{equation}
where $t$ is treated as a binomial element, that is, $p_k(t)=t$ for all $k$.
The inverse  series $\sigma_t[L(X)]$ gives the characters of the Eulerian idempotents  \cite[Th. 3.7]{Hanl}.

Otherwise said,
\begin{equation}\label{eq:ppol2}
	\sum_{i\ge 0}(-t)^i\ch H^{i}(P_n;\C) = \left. \lambda_{-1/t}[L(-tX)]\right|_{{\rm degree}\ n}.
\end{equation}
and we can extract a factor $\lambda_{-1/t}[\ell_1(-tX)]=\sigma_1(X)$. Each factor $\lambda_{-1/t}[\ell_k(-tX)]$
contains only positive powers of $t$, so that the coefficient of $t^i$ is a symmetric function of finite degree.
This proves the representation stability of $H^i(P_n,\C)$ in the sense of \cite{CF,CEF}. The calculation
of $H^2(P_\bullet,\C)$ presented 
in \cite[Example 5.1.A]{CEF} can be done as follows. The characteristic of $H^2(P_n,\C)$ is the coefficient of
$t^2$ in
\begin{equation}
	\begin{split}
		\lambda_{-1/t}[L(-tX)]&= \sigma_1\cdot\lambda_{-1/t}[\ell_2(-tX)+\ell_3(-tX)+\cdots]\\
		&=\sigma_1\left(1-\frac{1}{t}(\ell_2(-tX)+\ell_3(-tX))+\frac{1}{t^2}e_2[\ell_2(-tX)]+\cdots\right)\\
	\end{split}
\end{equation}
with $\ell_2(-X)=s_2$, $\ell_3(-X)=s_{21}$, $e_2[\ell_2(-X)]=s_{31}$
so that the coefficient of $t^2$ is
\begin{equation}
	\sigma_1\cdot(s_{21}+s_{31})=\llangle s_{21}+s_{31}\rrangle
\end{equation}
which is the character of the FI-module $M(21)+M(31)$ in the notation of \cite{CEF}.

For example,
\begin{align}
	\ch_t H^*(P_2)&=(1+t)s_2\\
	\ch_t H^*(P_3) &= s_3+t(s_3+s_{21})+t^2s_{21}\\
	\ch_t H^*(P_4) &= s_4+t(s_4+s_{31}+s_{22})+t^2(2s_{31}+s_{22}+s_{211})+t^3(s_{31}+s_{211})
\end{align}
and the coefficient of $t^2$ in the last equation is indeed the term of degree $4$ in  $\sigma_1(s_{21}+s_{31})$.

The pure braid group $P_n$ is the fundamental group of the variety
\begin{equation}
M_n=\{(z_1,\ldots,z_n)\in {\C}^n \ \mid\ z_i\not = z_j\ {\rm for}\
i\not = j \} \ .
\end{equation}
Arnold has shown that the cohomology $H^*(M_n,{\C})$ of this space is generated by the classes
$a_{ij}=[\omega_{ij}]$ of the holomorphic forms
\begin{equation}
\omega_{ij}={1\over 2\pi i}{dz_i-dz_j\over z_i-z_j}
\end{equation}
and is therefore isomorphic to the graded algebra
$A(n)$ generated over ${\C}$ by the elements
$a_{ij}=a_{ji}$ $i\not = j$ subject to the  relations
\begin{align}
&a_{ij}a_{rs}=-a_{rs}a_{ij}\\
&a_{ij}a_{jk}+a_{jk}a_{ki}+a_{ki}a_{ij}=0 \ . 
\end{align}
This is the so-called Arnold algebra, now a special case of an Orlik-Solomon algebra \cite{OrS}.

The natural action of $\SG_n$ on $M_n$, defined by
\begin{equation}
\sigma (z_1,\ldots,z_n) = (z_{\sigma(1)},\ldots,z_{\sigma(n)})
\end{equation}
induces an action of its cohomology, given by
$\sigma a_{ij}=a_{\sigma(i)\sigma(j)}$. The characteristic of this action
has been computed by Lehrer and Solomon \cite{LeS}, and their result is equivalent to Equation \eqref{eq:ppol}.

\section{Tables}
\subsection{Stable inner plethysms $\llb f\rrb = \hat f[\sigma_1h_1]$ in terms of stable permutation characters} 

\begin{align*}
\llb h_1\rrb  &= \llangle h_1 \rrangle \\
\llb h_2\rrb  &= \llangle h_2+h1  \rrangle\\
\llb h_{11}\rrb &= \llangle h_{11}+h_1 \rrangle\\
\llb h_3\rrb  &=  \llangle h_{3}+ h_{11}+ h_{1} \rrangle\\
\llb h_{21}\rrb  &=  \llangle  h_{21}+2 h_{11}+ h_{1} \rrangle\\
\llb h_{111}\rrb  &=  \llangle h_{111}+3 h_{11} h_{1}  \rrangle\\
\llb h_4\rrb  &=  \llangle  h_{4}+ h_{21}+ h_{2}+ h_{11}+ h_{1} \rrangle\\
\llb h_{31}\rrb  &=  \llangle  h_{31}+ h_{21}+ h_{111}+3 h_{11} +h_{1} \rrangle\\
\llb h_{22}\rrb  &=  \llangle  h_{22}+2 h_{21} +h_{111}+ h_{2}+3 h_{11}+ h_{1} \rrangle\\
\llb h_{211}\rrb  &=  \llangle  h_{211}+ h_{21}+3 h_{111}+5 h_{11}+h_{1} \rrangle\\
\llb h_{1111}\rrb  &=  \llangle  h_{1111}+6 h_{111}+7 h_{11} +h_{1} \rrangle
\end{align*}

\subsection{Stable permutation characters in terms of stable inner plethysms}
\begin{align*}
 \llangle h_{1} \rrangle &= \llb h_1 \rrb \\
 \llangle h_{2} \rrangle &= \llb h_{2} -h_{1} \rrb \\
 \llangle h_{11} \rrangle &= \llb  h_{11} -h_{1} \rrb \\
 \llangle h_{3} \rrangle &= \llb  h_{3}- h_{11}  \rrb \\
 \llangle h_{21} \rrangle &= \llb  h_{21}-2 h_{11}+ h_{1} \rrb \\
 \llangle h_{111} \rrangle &= \llb  h_{111}-3 h_{11}+2 h_{1} \rrb \\
 \llangle h_{4} \rrangle &= \llb  h_{4}- h_{21}+ h_{11}- h_{2} \rrb \\
 \llangle h_{31} \rrangle &= \llb  h_{31}- h_{21}- h_{111}+2 h_{11}- h_{1} \rrb \\
 \llangle h_{22} \rrangle &= \llb  h_{22}- 2h_{21}- h_{111}+4 h_{11}- h_{2}-h_{1}  \rrb \\
 \llangle h_{211} \rrangle &= \llb  h_{211} -h_{21}-3 h_{111}+6 h_{11}-3 h_{1} \rrb \\
 \llangle h_{1111} \rrangle &= \llb  h_{1111}-6 h_{111}+11 h_{11}-6 h_{1}\rrb  
\end{align*}

\subsection{Dual basis of $\tilde s_\lambda$, up to degree 5}
\begin{align*}
	\tilde s_1^* &= s_{1}+ s_{1 1} + 2s_{2}+ 3s_{2 1} + 4s_{3}+ s_{2 1 1} + 3s_{2 2} + 7s_{3 1} + 7s_{4} + 3s_{2 2 1} + 3s_{3 1 1} + 10s_{3 2} + 14s_{4 1} + 12s_{5}\\
	\tilde s_2^* &= s_{2}+ 2s_{2 1} + 2s_{3}+ s_{2 1 1} + 4s_{2 2} + 5s_{3 1} + 5s_{4}+ 4s_{2 2 1} + 4s_{3 1 1} + 11s_{3 2} + 13s_{4 1} + 9s_{5}\\
	\tilde s_{11}^*&= s_{1 1}+ s_{1 1 1} + 2s_{2 1} + s_{3}+ 3s_{2 1 1} + s_{2 2} + 6s_{3 1} + 2s_{4}+ s_{2 1 1 1} + 3s_{2 2 1} + 8s_{3 1 1} + 8s_{3 2} + 12s_{4 1} + 5s_{5}\\
	\tilde s_3^*&=  s_{3}+ s_{2 2} + 2s_{3 1} + 2s_{4}+ 2s_{2 2 1} + s_{3 1 1} + 6s_{3 2} + 6s_{4 1} + 5s_{5}\\
	\tilde s_{21}^*&=s_{2 1}+ 2s_{2 1 1} + 2s_{2 2} + 3s_{3 1} + s_{4}+ s_{2 1 1 1} + 5s_{2 2 1} + 7s_{3 1 1} + 8s_{3 2} + 9s_{4 1} + 3s_{5}\\
	\tilde s_{111}^*&=s_{1 1 1}+ s_{1 1 1 1} + 2s_{2 1 1} + s_{3 1}+ 3s_{2 1 1 1} + s_{2 2 1} + 6s_{3 1 1} + 2s_{3 2} + 3s_{4 1}\\
	\tilde s_4^*&= s_{4}+ s_{3 2} + 2s_{4 1} + 2s_{5}\\
	\tilde s_{31}^*&=s_{3 1}+ s_{2 2 1} + 2s_{3 1 1} + 3s_{3 2} + 3s_{4 1} + s_{5}\\
	\tilde s_{22}^*&=s_{2 2}+ 2s_{2 2 1} + s_{3 1 1} + 2s_{3 2} + s_{4 1}\\
	\tilde s_{211}^*&=s_{2 1 1}+ 2s_{2 1 1 1} + 2s_{2 2 1} + 3s_{3 1 1} + s_{3 2} + s_{4 1}\\
	\tilde s_{1111}^*&=s_{1 1 1 1}  s_{1 1 1 1 1} + 2s_{2 1 1 1} + s_{3 1 1}
\end{align*}

\subsection{Schur functions on the basis  $\tilde s_\mu$}
\begin{align*}
	s_{1} & =  \tilde s_{0} + \tilde s_{1}\\
s_{2} & =  2\tilde s_{0} + 2\tilde s_{1} + \tilde s_{2}\\
s_{1  1} & =  \tilde s_{1} + \tilde s_{1 1}\\
s_{3} & =  3\tilde s_{0} + 4\tilde s_{1} + \tilde s_{1 1} + 2\tilde s_{2} + \tilde s_{3}\\
s_{2  1} & =  \tilde s_{0} + 3\tilde s_{1} + 2\tilde s_{1 1} + 2\tilde s_{2} + \tilde s_{2 1}\\
s_{1  1  1} & =  \tilde s_{1 1} + \tilde s_{1 1 1}\\
s_{4} & =  5\tilde s_{0} + 7\tilde s_{1} + 2\tilde s_{1 1} + 5\tilde s_{2} + \tilde s_{2 1} + 2\tilde s_{3} + \tilde s_{4}\\
s_{3  1} & =  2\tilde s_{0} + 7\tilde s_{1} + 6\tilde s_{1 1} + \tilde s_{1 1 1} + 5\tilde s_{2} + 3\tilde s_{2 1} + 2\tilde s_{3} + \tilde s_{3 1}\\
s_{2  2} & =  2\tilde s_{0} + 3\tilde s_{1} + \tilde s_{1 1} + 4\tilde s_{2} + 2\tilde s_{2 1} + \tilde s_{2 2} + \tilde s_{3}\\
s_{2  1  1} & =  \tilde s_{1} + 3\tilde s_{1 1} + 2\tilde s_{1 1 1} + \tilde s_{2} + 2\tilde s_{2 1} + \tilde s_{2 1 1}\\
s_{1  1  1  1} & =  \tilde s_{1 1 1} + \tilde s_{1 1 1 1}
\end{align*}
\subsection{Dual basis of $\tilde h_\mu$ up to degree 5}

\begin{align*}
\tilde h_{1}^*  &=    m_{1}+ m_{11} + m_{2}+ m_{111} + m_{21} + m_{3}+ m_{1111} + m_{211} + m_{22} + m_{31} + m_{4}\\
	&+ m_{11111} + m_{2111} + m_{221} + m_{311} + m_{32} + m_{41} + m_{5}\\ 
\tilde h_{2}^*  &=    m_{2}+  m_{22} + m_{4}\\
\tilde h_{1 1}^*  &=   m_{11}+ 3m_{111} + 2m_{21} + m_{3}+ 7m_{1111} + 5m_{211} + 3m_{22} + 3m_{31} + m_{4}\\
	&+ 10m_{11111} + 7m_{2111} + 5m_{221} + 4m_{311} + 3m_{32} + 2m_{41} + m_{5}\\
\tilde h_{3}^*  &=    m_{3}\\
\tilde h_{2 1}^*  &=    m_{21}+ m_{211} + 2m_{22} + m_{31} + m_{4}+ m_{221} + m_{32} + m_{41} + m_{5}\\
\tilde h_{1 1 1}^*  &=   m_{111} 6m_{1111} + 3m_{211} + m_{22} + m_{31} 15m_{11111} + 9m_{2111} + 5m_{221} + 4m_{311} + 2m_{32} + m_{41}\\
\tilde h_{4}^*  &=    m_{4} \\
\tilde h_{3 1}^*  &=    m_{31}+ m_{311} + m_{32} + m_{41} + m_{5}\\
\tilde h_{2 2}^*  &=    m_{22}\\ 
\tilde h_{2 1 1}^*  &=   m_{211}+ 3m_{2111} + 4m_{221} + 2m_{311} + 2m_{32} + m_{41}\\
\tilde h_{1 1 1 1}^*  &=   m_{1111}+ 10m_{11111} + 4m_{2111} + m_{221} + m_{311} 
\end{align*}

\subsection{$h$ on $\tilde h$}

\begin{align*}
h_{1}  &=  \tilde h_{1}\\
h_{2}  &=  \tilde h_{1} + \tilde h_{2}\\
h_{1 1}  &=  \tilde h_{1} + \tilde h_{11}\\
h_{3}  &=  \tilde h_{1} + \tilde h_{11} + \tilde h_{3}\\
h_{2 1}  &=  \tilde h_{1} + 2\tilde h_{11} + \tilde h_{21}\\
h_{1 1 1}  &=  \tilde h_{1} + 3\tilde h_{11} + \tilde h_{111}\\
h_{4}  &=  \tilde h_{1} + \tilde h_{11} + \tilde h_{2} + \tilde h_{21} + \tilde h_{4}\\
h_{3 1}  &=  \tilde h_{1} + 3\tilde h_{11} + \tilde h_{111} + \tilde h_{21} + \tilde h_{31}\\
h_{2 2}  &=  \tilde h_{1} + 3\tilde h_{11} + \tilde h_{111} + \tilde h_{2} + 2\tilde h_{21} + \tilde h_{22}\\
h_{2 1 1}  &=  \tilde h_{1} + 5\tilde h_{11} + 3\tilde h_{111} + \tilde h_{21} + \tilde h_{211}\\
h_{1 1 1 1}  &=  \tilde h_{1} + 7\tilde h_{11} + 6\tilde h_{111} + \tilde h_{1111}
\end{align*}

\newpage
\footnotesize

\end{document}